# ASYMPTOTIC OF THE HEAT KERNEL IN GENERAL BENEDICKS DOMAINS

by

Pierre Collet °  Servet Martínez*, Jaime San Martín *

°C.N.R.S. UMR 7644, Physique Théorique, Ecole Polytechnique,
91128 Palaiseau Cedex, France.
e-mail: collet@cpht.polytechnique.fr
*Universidad de Chile, Facultad de Ciencias Físicas y Matemáticas,
Departamento de Ingeniería Matemática,
and Centro de Modelamiento Matemático CNRS UMR 2071,
Casilla 170-3 Correo 3, Santiago, Chile.
e-mail: smartine@dim.uchile.cl  jsanmart@dim.uchile.cl

*Abstract.* Using a new inequality relating the heat kernel and the probability of survival, we prove asymptotic ratio limit theorems for the heat kernel (and survival probability) in general Benedicks domains. In particular, the dimension of the cone of positive harmonic measures with Dirichlet boundary condition can be derived from the rate of convergence to zero of the heat kernel (or the survival probability).

*Keywords:* positive harmonic functions, ratio limit theorems, survival probability.

*AMS classification:* 60J65, 31B05

**I.Introduction and main results.**

There exits an increasing literature devoted to the study of the killed Brownian motion on unbounded domains and the associated survival probability, for instance see [C.Z.], [K.], [L.], [P.]. Here, we consider the large time asymptotic of the heat kernel in general Benedicks domains with Dirichlet boundary conditions. The case of Benedicks domains with compact holes was treated in [C.M.SM.1]. Although some cases with non compact holes can be studied using the same methods, it seems hopeless to extend this approach to the general case. In the present paper we derive a new inequality (of Nash type) which allows to control the general case of Benedicks domains. For these domains we answer positively a conjecture of Davies [D1.] about limits of ratios. The case where there is only one positive harmonic function vanishing on the boundary (up to a positive factor) was treated before for more general domains in [C.M.SM.1] and for exterior domains in [C.M.SM.2]. See also [B.J.]. The most interesting case is of course the case when the cone of non negative harmonic functions with Dirichlet boundary condition is of dimension larger than one. We will see that for Benedicks domains, the asymptotic rate of convergence to zero of the heat kernel (or of the survival probability) clearly reflects this property.

We recall that a Benedicks domain $\Omega$ is the complement of a closed set $D^c \times \{0\}$ where $D$ is an open non trivial subset of the hyperplane $x_d = 0$ (see [B.]). We also assume as in [B.] that each point of $D^c \times \{0\}$ is regular for the Dirichlet problem in $\Omega$. If $x \in \mathbf{R}^d$ we will sometimes write $x = (\vec{x}, x_d)$ with $\vec{x} \in \mathbf{R}^{d-1}$. We recall that for Benedicks domains the cone $\mathcal{P}_\Omega$ of non negative harmonic functions with zero boundary condition is of dimension one or two (see [B.]). An integral test for this dimension is given in [B.] (see also [A.]) while Theorem 1 below gives a probabilistic characterization. In the case of dimension two, we denote by $u_1$ the extremal of the cone normalized by

$$\lim_{x_d \to +\infty} \frac{u_1(\vec{x}, x_d)}{x_d} = 1 \ .$$

We recall that since $u_1$ is an extremal we have

$$\lim_{x_d \to -\infty} \frac{u_1(\vec{x}, x_d)}{x_d} = 0 \ .$$

We denote by $u_2$ the symmetric function with respect to the hyperplane $x_d = 0$, namely

$$u_2(\vec{x}, x_d) = u_1(\vec{x}, -x_d).$$

It is another extremal function which satisfies

$$\lim_{x_d \to +\infty} \frac{u_2(\vec{x}, x_d)}{x_d} = 0 \quad \text{and} \quad \lim_{x_d \to -\infty} \frac{u_2(\vec{x}, x_d)}{x_d} = 1 \ .$$

Note that $u_1$ and $u_2$ generate $\mathcal{P}_\Omega$. We also denote by $v_s$ the function

$$v_s(\vec{x}, x_d) = u_1(\vec{x}, x_d) + u_2(\vec{x}, x_d) \ ,$$

which is modulo a constant factor the unique symmetric element in $\mathcal{P}_\Omega$.

We denote by $p_t(x, y)$ the heat kernel of $\Omega$ with Dirichlet boundary conditions. For a Brownian motion $X_t$ in $\Omega$ with absorption at the boundary, let $T_\Omega$ denote the exit time from $\Omega$. We will denote this random variable by $T$ when there is no ambiguity on the domain. $\mathbf{P}_x(T > t)$ will denote the probability of survival up to time $t$ starting at the point $x \in \Omega$, namely

$$\mathbf{P}_x(T > t) = \int_\Omega p_t(x, y) dy \; .$$

We will also use the heat kernel of the half space $x_d > 0$ denoted by $p_t^H(x, y)$. We now formulate our result.

**Theorem 1.** *For any Benedicks domain $\Omega$, for each $(x, y) \in \Omega \times \Omega$, the function of time $t^{1+d/2} p_t(x, y)$ has a limit when $t$ tends to infinity.*
*If this limit is infinite for one pair $(x, y)$, it is infinite for any pair $(x', y') \in \Omega \times \Omega$, the cone $\mathcal{P}_\Omega$ is of dimension one, and we have for any $x_1, \cdots, x_4$ in $\Omega$*

$$\lim_{t \to \infty} \frac{p_t(x_1, x_2)}{p_t(x_3, x_4)} = \frac{u(x_1) u(x_2)}{u(x_3) u(x_4)}$$

*where $u$ is any function in the cone.*
*If the limit is finite for one pair $(x, y)$ (and hence for all pairs of points), the cone $\mathcal{P}_\Omega$ has dimension two and we have for any $x$ and $y$ in $\Omega$*

$$\lim_{t \to \infty} t^{1+d/2} p_t(x, y) = \frac{2}{(2\pi)^{d/2}} \left( u_1(x) u_1(y) + u_2(x) u_2(y) \right) \; .$$

In the first case, the exact asymptotic rate in time of $p_t(x, y)$ differs according to the domain. For example, in dimension two, for the exterior of a segment one gets $t^{-1}(\log t)^{-2}$, and for the complement of a half line one gets $t^{-3/2}$ (see [C.M.SM.2] and [B.S.]).

A result similar to Theorem 1 holds for the survival probability.

**Theorem 2.** *For any Benedicks domain $\Omega$, for each $x \in \Omega$, the function of time $t^{1/2} \mathbf{P}_x(T > t)$ has a limit when $t$ tends to infinity.*
*If this limit is infinite for one point $x \in \Omega$, it is infinite for any point $y \in \Omega$. The cone $\mathcal{P}_\Omega$ is of dimension one, and we have for any $x$, and $y$ in $\Omega$*

$$\lim_{t \to \infty} \frac{\mathbf{P}_x(T > t)}{\mathbf{P}_y(T > t)} = \frac{u(x)}{u(y)}$$

*where $u$ is any function in the cone.*
*If the limit is finite (for one point), the cone $\mathcal{P}_\Omega$ has dimension two and we have for any $x \in \Omega$*

$$\lim_{t \to \infty} t^{1/2} \mathbf{P}_x(T > t) = \sqrt{\frac{2}{\pi}} v_s(x) \; .$$

Here also, in the first case the asymptotic rate depends on the domain.

The rest of the paper is organized as follows. In section II we prove Theorem 1. The proof of Theorem 2 is given in section III. In the appendix we give a proof of symmetry estimates needed in the proof of Lemma 4.

## II. Asymptotic of the heat kernel.

The proof of Theorem 1 will be a consequence of several results. We first establish a general estimate which may be of independent interest. This Lemma follows by applying Nash inequality to the heat kernel. In fact, an analytic proof is based on the inequality $\|\phi\|_2^4 \leq C_2 \|\nabla \phi\|_2^2 \|\phi\|_1^2$, for $\phi(x) = p_t(x, y)$. A direct computation yields $\|p_t(\cdot, y)\|_2^2 = p_{2t}(y, y)$, $\|\nabla . p_t(\cdot, y)\|_2^2 = -\dot{p}_{2t}(y, y)$ and $\|p_t(\cdot, y)\|_1 = \mathbf{P}_y(T > t)$, leading to an inequality of the type (1) below, with $p_{4t}$ instead of $p_{3t}$. In the sequel we give a direct and simpler probabilistic proof.

**Lemma 3.** *For any $x$ and $y$ in $\Omega$ and any $t > 0$, we have*

$$p_{3t}(x, y) \leq \frac{1}{(2\pi t)^{d/2}} \mathbf{P}_x(T > t) \mathbf{P}_y(T > t) . \tag{1}$$

**Proof**. By the semi group property we have

$$p_{3t}(x, y) = \int_\Omega dz_1 \int_\Omega dz_2 \; p_t(x, z_1) p_t(z_1, z_2) p_t(z_2, y) .$$

On the other hand from the Gaussian bound

$$p_t(x, y) \leq \frac{1}{(2\pi t)^{d/2}} e^{-|x-y|^2/2t}$$

we obtain

$$p_t(z_1, z_2) \leq \frac{1}{(2\pi t)^{d/2}} e^{-|z_1-z_2|^2/2t} \leq \frac{1}{(2\pi t)^{d/2}} .$$

The result follows from the above expression for $\mathbf{P}_x(T > t)$. QED

By a direct computation, using the results of [B.S.] or [D2.], one can check that for cones the following limit exists

$$\lim_{t \to \infty} \frac{(2\pi t)^{d/2} p_t(x_1, x_2)}{\mathbf{P}_{x_1}(T > t) \mathbf{P}_{x_2}(T > t)} .$$

This limit also exists in the case of exterior domains as follows from the results in [C.M.SM.2]. For the case of Benedicks domains with two extremal harmonic functions this follows from Theorems 1 and 2.

**Lemma 4.** *For any Benedicks domain $\Omega$ with a cone $\mathcal{P}_\Omega$ of dimension two, there is a constant $\Gamma > 1$ such that for any $t > 1$ and any $(x, y) \in \Omega \times \Omega$ we have*

$$p_t(x, y) \leq \frac{\Gamma}{t^{1+d/2}} \big(1 + v_s(x)\big)\big(1 + v_s(y)\big) \tag{2}$$

*We also have for some finite constant $K$, for any $t > 0$ and any $x \in \Omega$ the bound*

$$\mathbf{P}_x(T > t) \leq K \frac{1 + v_s(x)}{\sqrt{t}} .$$

**Proof.** We first observe that if $\Gamma > 4/(2\pi)^{d/2}$, the first estimate is true for $1 \leq t \leq 4$ due to the Gaussian bound and the non negativity of $v_s$. We are now going to prove recursively that if (2) holds on the time interval $[1, 3^n]$ (where $n$ is a positive integer) it also holds with the same constant on the time interval $[3^n, 3^{n+1}]$.

We will assume that $y_d > 0$, the case $y_d < 0$ is similar, and $y_d = 0$ is recovered by continuity. We first prove an upper bound for $\mathbf{P}_y(T > t)$, using the recursive bound (2) (i.e. assuming $t \in [1, 3^n]$). Let $\alpha = \Gamma^{-1/(2+d)}$. We have

$$\mathbf{P}_y(T > t) = I_1 + I_2 \;,$$

where

$$I_1 = \int_{|z_d| \leq \alpha \sqrt{t}} p_t(z, y) dz \quad \text{and} \quad I_2 = \int_{|z_d| \geq \alpha \sqrt{t}} p_t(z, y) dz \;.$$

We first estimate $I_2$. Recall that (see [B.]) $v_s(z) \geq |z_d|$. Therefore

$$I_2 = \int_{|z_d| \geq \alpha \sqrt{t}} \frac{1 + v_s(z)}{1 + v_s(z)} p_t(z, y) dz \leq \frac{1}{1 + \alpha \sqrt{t}} \int_{|z_d| \geq \alpha \sqrt{t}} (1 + v_s(z)) p_t(z, y) dz$$

$$\leq \frac{1}{1 + \alpha \sqrt{t}} \int (1 + v_s(z)) p_t(z, y) dz \leq \frac{1 + v_s(y)}{\alpha \sqrt{t}} \;.$$

In order to estimate $I_1$ we will split this integral into two parts. Namely

$$I_1 = J_1 + J_2$$

where

$$J_1 = \int_{|z_d| \leq \alpha \sqrt{t}, |\bar{z} - \bar{y}| \leq \alpha \sqrt{d-1} \sqrt{t}} p_t(z, y) dz \quad \text{and} \quad J_2 = \int_{|z_d| \leq \alpha \sqrt{t}, |\bar{z} - \bar{y}| \geq \alpha \sqrt{d-1} \sqrt{t}} p_t(z, y) dz \;.$$

We first estimate $J_2$. We observe that

$$J_2 \leq \sum_{l=1}^{d-1} J_{2,l} \quad \text{with} \quad J_{2,l} = \int_{|z_d| \leq \alpha \sqrt{t}, |z_l - y_l| \geq \alpha \sqrt{t}} p_t(z, y) dz \;.$$

We now estimate each $J_{2,l}$ in terms of $I_2$. To this end we use the symmetry with respect to the hyperplanes $z_d = \pm(z_l - y_l)$ denoted by $S_l^\pm$. Using Lemma A from the appendix we get for $z_l - y_l < -|z_d|$

$$p_t(z, y) \leq p_t(S_l^+ z, y) + p_t(S_l^- z, y)$$

and a similar estimate if $z_l - y_l > |z_d|$. Therefore

$$J_{2,l} \leq 4 I_2$$

and we get a bound by the previous estimate of $I_2$.

We now estimate $J_1$. We have

$$J_1 = \int_{|z_d|\leq\alpha\sqrt{t}\,,\,|\vec{z}-\vec{y}|\leq\alpha\sqrt{d-1}\sqrt{t}} p_t(z,y)dz =$$

$$\int_{|z_d|\leq\alpha\sqrt{t}\,,\,|\vec{z}-\vec{y}|\leq\alpha\sqrt{d-1}\sqrt{t}} p_t(z,y)^{1/2}(1+v_s(z))^{1/2}\frac{p_t(z,y)^{1/2}}{(1+v_s(z))^{1/2}}dz\,,$$

and using Schwartz inequality we get

$$J_1 \leq$$

$$\left(\int_{|z_d|\leq\alpha\sqrt{t}\,,\,|\vec{z}-\vec{y}|\leq\alpha\sqrt{d-1}\sqrt{t}} p_t(z,y)(1+v_s(z))dz \int_{|z_d|\leq\alpha\sqrt{t}\,,\,|\vec{z}-\vec{y}|\leq\alpha\sqrt{d-1}\sqrt{t}} \frac{p_t(z,y)}{(1+v_s(z))}dz\right)^{1/2}$$

$$\leq (1+v_s(y))^{1/2}\left(\int_{|z_d|\leq\alpha\sqrt{t}\,,\,|\vec{z}-\vec{y}|\leq\alpha\sqrt{d-1}\sqrt{t}} \frac{p_t(z,y)}{(1+v_s(z))}dz\right)^{1/2}.$$

In the last integral we use the recursive assumption (2) and get

$$J_1 \leq (1+v_s(y))\Gamma^{1/2}t^{-(1/2+d/4)}\left(\int_{|z_d|\leq\alpha\sqrt{t}\,,\,|\vec{z}-\vec{y}|\leq\alpha\sqrt{d-1}\sqrt{t}} \frac{1+v_s(z)}{(1+v_s(z))}dz\right)^{1/2}$$

$$\leq C_d(1+v_s(y))\Gamma^{1/2}\alpha^{d/2}t^{-1/2}\,,$$

where $C_d$ is a constant depending only on the dimension $d$.

Collecting the estimates and using the relation $\alpha = \Gamma^{-1/(2+d)}$, we get for a positive constant $C_4$ (which can always be assumed larger than 1) independent of $\Gamma$, $x$, $y$ and $t \in [1, 3^n]$

$$\mathbf{P}_y(T > t) \leq \frac{C_4(1+v_s(y))}{\sqrt{t}}\left(\frac{1}{\alpha}+\alpha^{d/2}\Gamma^{1/2}\right) \leq 2\frac{C_4(1+v_s(y))}{\sqrt{t}}\Gamma^{1/(2+d)}\,.$$

This is the second part of the Lemma for $t \in [1, 3^n]$ if we take

$$K = 2C_4\Gamma^{1/(2+d)}\,.$$

For $t \in [0,1]$ the second part of the Lemma follows at once if we take $K > 1$ since $\mathbf{P}_y(T > t) \leq 1$. Using Lemma 3 we get for $t \in [3^n, 3^{n+1}]$

$$p_t(x,y) \leq \frac{4}{(2\pi)^{d/2}}C_4^2\Gamma^{2/(2+d)}(1+v_s(x))(1+v_s(y))\frac{1}{t^{1+d/2}}\,.$$

If we now chose

$$\Gamma > \frac{(2C_4)^{(4+2d)/d}}{(2\pi)^{(2+d)/2}}$$

we obtain for $t \in [3^n, 3^{n+1}]$ the required estimate, and we complete recursively the proof of the first part of Lemma 4. QED

We now come to the proof of Theorem 1. First of all, we observe that if there is a pair $(x, y) \in \Omega \times \Omega$ such that $t^{(1+d/2)} p_t(x, y)$ has a finite accumulation point when $t$ tends to infinity, then using Harnack inequality we conclude that $t^{(1+d/2)} p_t(x, x)$ has also a finite accumulation point. This implies as in [C.M.SM.1] that we can find a diverging sequence $(t_n)$ such that for any $s > 0$ and any $(y, z) \in \Omega \times \Omega$, the sequence $(t_n + s)^{(1+d/2)} p_{t_n+s}(y, z)$ converges to a function $\varphi(y, z)$ which is harmonic in $y$ and $z$ (separately) and vanishes at the boundary of $\Omega$. In the case of manifolds, see [A.B.J.] for a different proof. Note that at this point of the argument, $\varphi$ may still depend on the sequence $(t_n)$. We will prove below that this is not the case, and this implies convergence. We recall that the explicit formula for $p_t^H(x, y)$ is

$$p_t^H(x, y) = \frac{2}{(2\pi t)^{d/2}} e^{-|\vec{x}-\vec{y}|^2/2t} e^{-(x_d^2+y_d^2)/2t} \sinh\left(\frac{x_d y_d}{t}\right) . \tag{3}$$

It follows immediately that

$$\lim_{t \to \infty} t^{1+d/2} p_t^H(x, y) = \frac{2 x_d y_d}{(2\pi)^{d/2}} . \tag{4}$$

From the bound

$$p_t^H(x, y) \leq p_t(x, y)$$

we conclude that

$$\varphi(y, z) \geq \frac{2 y_d z_d}{(2\pi)^{d/2}} .$$

In particular, for fixed $z$ we find that $\varphi(y, z)$ is a non negative harmonic function vanishing on the boundary and growing linearly with $y_d$. This implies by the results of [B.] that the cone $\mathcal{P}_\Omega$ is of dimension two.

Using the various symmetries and the convexity of the cone $\mathcal{P}_\Omega$ we can say more about the structure of the function $\varphi(x, y)$. First of all since $\mathcal{P}_\Omega$ is of dimension two with $u_1$ and $u_2$ two extremals, and since $\varphi(x, y)$ is harmonic in $x$ at fixed $y$, non negative and zero on the boundary, we can find two functions $f_1$ and $f_2$ of $y$ such that

$$\varphi(x, y) = f_1(y) u_1(x) + f_2(y) u_2(x) .$$

Since $u_1$ and $u_2$ are not proportional, we conclude that $f_1(y)$ and $f_2(y)$ are harmonic, non negative and zero on the boundary of $\Omega$. Therefore writing $f_1$ and $f_2$ as nonnegative combinations of $u_1$ and $u_2$, and using also the symmetry of $\varphi$ in $x$ and $y$ we obtain

$$\varphi(x, y) = a u_1(x) u_1(y) + b u_2(x) u_2(y) + c\big(u_1(x) u_2(y) + u_2(x) u_1(y)\big) ,$$

where $a$, $b$ and $c$ are three nonnegative constants which may depend on the sequence $(t_n)$. Using the symmetry with respect to the hyperplane $x_d = 0$, we conclude that $a = b$.

We now use the formula (see [C.M.SM.1]) for $x_d y_d > 0$

$$p_t(x,y) = p_t^H(x,y) - \frac{1}{2}\int_0^t ds \int_D \partial_{n_{\vec{\xi}}} p_{t-s}^H(x,(\vec{\xi},0)) p_s((\vec{\xi},0),y) d\vec{\xi}$$

$$= p_t^H(x,y) + \frac{x_d}{(2\pi)^{d/2}} \int_0^t \frac{ds}{(t-s)^{1+d/2}} \int_D e^{-x_d^2/2(t-s)} e^{-|\vec{x}-\vec{\xi}|^2/2(t-s)} p_s((\vec{\xi},0),y) d\vec{\xi}. \quad (5)$$

If $x_d y_d < 0$ the same identity without the term $p_t^H(x,y)$ holds. Notice that the function $p_s((\xi,0),y)$ is symmetric with respect to the hyperplane $y_d = 0$. Therefore if $y^*$ is the symmetric of $y$ with respect to this hyperplane, we have for $y_d x_d > 0$

$$p_t(x,y) = p_t^H(x,y) + p_t(x,y^*) .$$

Since $u_1(y^*) = u_2(y)$ and $u_2(y^*) = u_1(y)$ we conclude from (4) and our normalization of $u_1$ (and $u_2$) that $a = 2(2\pi)^{-d/2} - c$. It now remains to show that $c = 0$.

We now assume $x_d > 0$ and split the integral

$$x_d \int_0^t \frac{ds}{(t-s)^{1+d/2}} \int_D e^{-x_d^2/2(t-s)} e^{-|\vec{x}-\vec{\xi}|^2/2(t-s)} p_s((\vec{\xi},0),y) d\vec{\xi} = K_1 + K_2 ,$$

in two pieces given by

$$K_1 = x_d \int_0^{t/2} \frac{ds}{(t-s)^{1+d/2}} \int_D e^{-x_d^2/2(t-s)} e^{-|\vec{x}-\vec{\xi}|^2/2(t-s)} p_s((\vec{\xi},0),y) d\vec{\xi},$$

$$K_2 = x_d \int_{t/2}^t \frac{ds}{(t-s)^{1+d/2}} \int_D e^{-x_d^2/2(t-s)} e^{-|\vec{x}-\vec{\xi}|^2/2(t-s)} p_s((\vec{\xi},0),y) d\vec{\xi}.$$

We now estimate each term separately starting with $K_2$. Using Lemma 4 we have

$$K_2 \leq \Gamma(1+v_s(y)) \int_{t/2}^t \frac{ds}{s^{1+d/2}(t-s)^{1+d/2}} \int_D e^{-x_d^2/2(t-s)} (1+v_s(\vec{\xi},0)) e^{-|\vec{x}-\vec{\xi}|^2/2(t-s)} d\vec{\xi}$$

$$\leq \Gamma 2^{1+d/2} t^{-(1+d/2)}(1+v_s(y)) \int_{t/2}^t \frac{ds}{(t-s)^{1+d/2}} \int_D e^{-x_d^2/2(t-s)} (1+v_s(\vec{\xi},0)) e^{-|\vec{x}-\vec{\xi}|^2/2(t-s)} d\vec{\xi}.$$

We integrate first over $s$ and obtain

$$K_2 \leq \mathcal{O}(1) t^{-(1+d/2)}(1+v_s(y)) x_d \int_D \frac{1+v_s(\vec{\xi},0)}{(x_d^2 + |\vec{\xi}-\vec{x}|^2)^{d/2}} d\vec{\xi} .$$

Using the Herglotz representation for harmonic functions in a half hyperplane, we obtain finally

$$K_2 \leq C(\vec{x}) t^{-(1+d/2)}(1+v_s(y))(1+o(x_d)) ,$$

where $C(\vec{x})$ is bounded on compact sets in $\mathbf{R}^{d-1}$.

We now estimate $K_1$. We first observe that

$$t^{-(1+d/2)} x_d \int_0^{t/2} ds \int_D e^{-x_d^2/t} e^{-|\vec{x}-\vec{\xi}|^2/t} p_s((\vec{\xi}, 0), y) d\vec{\xi} \leq K_1 \leq (2\pi)^{d/2} p_t(x, y) .$$

Therefore, multiplying by $t^{1+d/2}$, taking $t = t_n$ and letting $n$ tend to infinity, we obtain by the monotone convergence Theorem

$$x_d \int_0^\infty ds \int_D p_s((\vec{\xi}, 0), y) d\vec{\xi}$$
$$\leq (2\pi)^{d/2} (a(u_1(x) u_1(y) + u_2(x) u_2(y)) + c(u_1(x) u_2(y) + u_2(x) u_1(y))) .$$

Dividing by $x_d$ and then taking $x_d$ to infinity, we obtain

$$V(y) = \int_0^\infty ds \int_D p_s((\vec{\xi}, 0), y) d\vec{\xi} \leq (2\pi)^{d/2} (a + c) v_s(y) .$$

In particular, this function is zero on the boundary of $\Omega$. We now derive a better estimate on $V$.

In equation (5) we take $y = (\vec{\eta}, 0)$ and we integrate over $\vec{\eta}$ and $t$. For any $\epsilon > 0$, we obtain

$$\int_\epsilon^{\epsilon^{-1}} dt \int_{D, |\vec{\eta}| \leq \epsilon^{-1}} p_t(x, (\vec{\eta}, 0)) d\vec{\eta}$$
$$= \frac{x_d}{(2\pi)^{d/2}} \int_\epsilon^{\epsilon^{-1}} dt \int_0^t \frac{ds}{(t-s)^{1+d/2}} \int_D d\vec{\xi} e^{-(|\vec{\xi}-\vec{x}|^2 + x_d^2)/2(t-s)} \int_{D, |\vec{\eta}| \leq \epsilon^{-1}} p_s((\vec{\xi}, 0), (\vec{\eta}, 0)) d\vec{\eta} .$$

We now interchange the integrations in $s$ and $t$ and obtain

$$\int_\epsilon^{\epsilon^{-1}} dt \int_{D, |\vec{\eta}| \leq \epsilon^{-1}} p_t(x, (\vec{\eta}, 0)) d\vec{\eta}$$
$$= \frac{x_d}{(2\pi)^{d/2}} \int_0^{\epsilon^{-1}} ds \int_{\max\{\epsilon, s\}}^{\epsilon^{-1}} \frac{dt}{(t-s)^{1+d/2}} \int_D d\vec{\xi} e^{-(|\vec{\xi}-\vec{x}|^2 + x_d^2)/2(t-s)} \int_{D, |\vec{\eta}| \leq \epsilon^{-1}} p_s((\vec{\xi}, 0), (\vec{\eta}, 0)) d\vec{\eta} ,$$
$$\leq \frac{x_d}{(2\pi)^{d/2}} \int_0^{\epsilon^{-1}} ds \int_0^\infty \frac{d\tau}{\tau^{1+d/2}} \int_D d\vec{\xi} e^{-(|\vec{\xi}-\vec{x}|^2 + x_d^2)/2\tau} \int_{D, |\vec{\eta}| \leq \epsilon^{-1}} p_s((\vec{\xi}, 0), (\vec{\eta}, 0)) d\vec{\eta}$$
$$\leq \mathcal{O}(1) x_d \int_D d\vec{\xi} \frac{1}{(|\vec{\xi} - \vec{x}|^2 + x_d^2)^{d/2}} \int_0^{\epsilon^{-1}} ds \int_{D, |\vec{\eta}| \leq \epsilon^{-1}} p_s((\vec{\xi}, 0), (\vec{\eta}, 0)) d\vec{\eta} .$$

Since $V(y)$ is a non negative function (in fact harmonic) in the upper half hyperplane, we can let $\epsilon$ tend to zero and using the monotone convergence Theorem we obtain

$$V(x) \leq \mathcal{O}(1) x_d \int_D d\vec{\xi} \frac{V(\vec{\xi}, 0)}{(|\vec{\xi} - \vec{x}|^2 + x_d^2)^{d/2}} \leq \mathcal{O}(1) x_d \int_D d\vec{\xi} \frac{v_s(\vec{\xi}, 0)}{(|\vec{\xi} - \vec{x}|^2 + x_d^2)^{d/2}} .$$

This is again the Herglotz representation and implies $V(x) = C(\vec{x})(1 + o(x_d))$, where $C(\vec{x})$ is bounded on compact sets in $\mathbf{R}^{d-1}$.

We now come back to the estimation of $K_1$. We have

$$K_1 \leq 2^{1+d/2} t^{-(1+d/2)} x_d \int_0^{t/2} ds \int_D p_s((\vec{\xi}, 0), y) d\vec{\xi}$$

$$\leq x_d t^{-(1+d/2)} \mathcal{O}(1) V(y) \leq x_d t^{-(1+d/2)} C(\vec{y})(1 + o(y_d)) \,.$$

We now take $y_d = -x_d < 0$ large, $\vec{x} = \vec{y}$ fixed. We have

$$\lim_{n \to \infty} t_n^{1+d/2} p_{t_n}(x, y) = \varphi(x, y)$$

$$= 2a u_1(x) u_2(x) + c(u_1(x)^2 + u_2(x)^2) = c\, x_d^2 + o(1 + x_d^2) \leq \mathcal{O}(1) x_d (1 + o(x_d)) \,,$$

which implies c=0. Therefore $a = 2(2\pi)^{-d/2}$ for any accumulation point, which implies convergence and finishes the proof of the last part of Theorem 1.

Finally, if $t^{(1+d/2)} p_t(x, y)$ diverges for a pair of points $(x, y)$ in $\Omega$, by Harnack inequality, it diverges for all pair of points. It follows from Lemma 4 that the cone $\mathcal{P}_\Omega$ is of dimension one and the first part of Theorem 1 follows as in [C.M.SM.1]. This finishes the proof of Theorem 1. QED

### III. Asymptotic of the survival probability.

In this section we give a proof of Theorem 2 on the asymptotic of the survival probability.

We first consider the case where for some $x \in \Omega$ the function of time $t^{1/2} \mathbf{P}_x(T > t)$ has a finite accumulation point. In other words, there is a diverging sequence $(t_n)$ such that $t_n^{1/2} \mathbf{P}_x(T > t_n)$ has a finite limit. Applying Lemma 3, we conclude that $t_n^{1+d/2} p_{t_n}(x, x)$ is bounded. Therefore we are in the second case of Theorem 1, namely the cone $\mathcal{P}_\Omega$ is of dimension 2. Applying Lemma 4 we conclude that $t^{1/2} \mathbf{P}_x(T > t)$ is bounded uniformly in $t$. On the other hand, any accumulation point is a non negative harmonic function which is symmetric with respect to the hyperplane $x_d = 0$ and zero on the boundary of $\Omega$. This last statement follows (for $t > 1$) from Lemma 4 by the estimate

$$\sqrt{t} \mathbf{P}_x(T > t) = \int_\Omega p_1(x, y) \sqrt{t} \mathbf{P}_y(T > t - 1)\, dy$$

$$\leq K \sqrt{\frac{t}{t-1}} \int_\Omega p_1(x, y)(1 + v_s(y))\, dy \leq K \sqrt{\frac{t}{t-1}} (\mathbf{P}_x(T \geq 1) + v_s(x)) \,.$$

Therefore, any accumulation point is of the form $A v_s(x)$ where the constant $A$ may depend on the accumulation point. We now show that this is not the case. Assuming $x_d > 0$ and integrating relation (5) with respect to $y \in \Omega$ we obtain the following expansion for $\mathbf{P}_x(T > t)$

$$\mathbf{P}_x(T_H > t) + \frac{x_d}{(2\pi)^{d/2}} \int_0^t \frac{ds}{(t-s)^{1+d/2}} \int_D e^{-(|\vec{x}-\vec{\xi}|^2 + x_d^2)/2(t-s)} \mathbf{P}_{(\vec{\xi}, 0)}(T > s) d\vec{\xi} \,.$$

We also have easily from (3)

$$\mathbf{P}_x(T_H > t) = \frac{1}{\sqrt{2\pi t}} \int_{-x_d}^{x_d} e^{-z^2/2t} dz = \frac{2x_d}{\sqrt{2\pi t}} + \mathcal{O}((x_d/\sqrt{t})^3) \, .$$

We now split the integral in the above expression of $\mathbf{P}_x(T > t)$ in two parts, namely

$$\frac{x_d}{(2\pi)^{d/2}} \int_0^t \frac{ds}{(t-s)^{1+d/2}} \int_D e^{-|\vec{x}-\vec{\xi}|^2/2(t-s)} e^{-x_d^2/2(t-s)} \mathbf{P}_{(\vec{\xi},0)}(T > s) d\vec{\xi} = M_1 + M_2$$

with

$$M_1 = \frac{x_d}{(2\pi)^{d/2}} \int_0^{t/2} \frac{ds}{(t-s)^{1+d/2}} \int_D e^{-|\vec{x}-\vec{\xi}|^2/2(t-s)} e^{-x_d^2/2(t-s)} \mathbf{P}_{(\vec{\xi},0)}(T > s) d\vec{\xi} \, ,$$

and

$$M_2 = \frac{x_d}{(2\pi)^{d/2}} \int_{t/2}^t \frac{ds}{(t-s)^{1+d/2}} \int_D e^{-|\vec{x}-\vec{\xi}|^2/2(t-s)} e^{-x_d^2/2(t-s)} \mathbf{P}_{(\vec{\xi},0)}(T > s) d\vec{\xi} \, .$$

We first estimate $M_1$. We have, using Lemma 4, that $M_1$ is bounded above by

$$\mathcal{O}(1) x_d \int_0^{t/2} \frac{ds}{\sqrt{s}(t-s)} \int_D e^{-(|\vec{x}-\vec{\xi}|^2+x_d^2)/2(t-s)} \left(\frac{|\vec{x}-\vec{\xi}|^2+x_d^2}{2(t-s)}\right)^{d/2} \frac{1+v_s(\vec{\xi},0)}{(|\vec{x}-\vec{\xi}|^2+x_d^2)^{d/2}} d\vec{\xi}$$

$$\leq \mathcal{O}(1) x_d \int_0^{t/2} \frac{ds}{\sqrt{s}(t-s)} \int_D \frac{1+v_s(\vec{\xi},0)}{(|\vec{x}-\vec{\xi}|^2+x_d^2)^{d/2}} d\vec{\xi} \leq \mathcal{O}(1) x_d t^{-1/2} \int_D \frac{1+v_s(\vec{\xi},0)}{(|\vec{x}-\vec{\xi}|^2+x_d^2)^{d/2}} \, ,$$

and since this last integral is proportional to the Herglotz integral, we obtain again

$$M_1 \leq t^{-1/2} C(\vec{x})(1 + o(x_d)) \, ,$$

where $C(\vec{x})$ is bounded on compact sets in $\mathbf{R}^{d-1}$. For $M_2$ we have using again Lemma 4

$$M_2 \leq \mathcal{O}(1) t^{-1/2} x_d \int_{t/2}^t \frac{ds}{(t-s)^{1+d/2}} \int_D e^{-(|\vec{x}-\vec{\xi}|^2+x_d^2)/2(t-s)} (1+v_s(\vec{\xi},0)) d\vec{\xi} \, .$$

Integrating first over $s$, we obtain again the Herglotz integral and conclude as before that

$$M_2 \leq \mathcal{O}(1) t^{-1/2} C(\vec{x})(1 + o(x_d)) \, .$$

Finally we obtain for fixed $\vec{x}$

$$\mathbf{P}_x(T > t) = \frac{2x_d t^{-1/2}}{\sqrt{2\pi}} + \mathcal{O}((x_d/\sqrt{t})^3) + \mathcal{O}(1) C(\vec{x}) t^{-1/2}(1 + o(x_d)) \, ,$$

which implies $A = \sqrt{2/\pi}$.

Assume now that $t^{1/2}\mathbf{P}_x(T > t)$ tends to infinity for some $x \in \Omega$. By Harnack inequality, this also holds for any $y \in \Omega$. Moreover from Lemma 4 we are in the case where the cone $\mathcal{P}_\Omega$ is of dimension one.

We now show that for any $s \in \mathbf{R}$ and any $x \in \Omega$, we have

$$\lim_{t \to \infty} \frac{\mathbf{P}_x(T > t + s)}{\mathbf{P}_x(T > t)} = 1 \ .$$

It is enough to consider $s > 0$ since the case $s < 0$ is obtained by taking the inverse of the above expression and replacing $t$ by $t + s$. By the Gaussian bound and integrating first over $\vec{y}$ we have

$$\int_{|y_d| \leq t^{1/4}\mathbf{P}_x(T>t)^{1/2}} p_t(x,y) dy \leq \frac{1}{(2\pi t)^{1/2}} \int_{|y_d| \leq t^{1/4}\mathbf{P}_x(T>t)^{1/2}} dy_d = \sqrt{\frac{2}{\pi}} \frac{\mathbf{P}_x(T > t)^{1/2}}{t^{1/4}} \ .$$

This implies

$$\mathbf{P}\left(|X_t^d| \leq t^{1/4}\mathbf{P}_x(T > t)^{1/2} \,\big|\, T > t\right) \leq \sqrt{\frac{2}{\pi}} \frac{1}{\left(\sqrt{t}\mathbf{P}_x(T > t)\right)^{1/2}} \ ,$$

which tends to zero when $t$ tends to infinity by our previous assumption.

For $s > 0$ the previous estimate implies

$$\frac{\mathbf{P}_x(T > t + s)}{\mathbf{P}_x(T > t)} = \mathbf{E}_x\left(\mathbf{P}_{X_t}(T > s) \,\big|\, T > t\right)$$

$$\geq \mathbf{E}_x\left(\mathbf{P}_{X_t}(T > s) \,\big|\, |X_t^d| \geq t^{1/4}\mathbf{P}_x(T > t)^{1/2}, T > t\right) \mathbf{P}_x\left(|X_t^d| \geq t^{1/4}\mathbf{P}_x(T > t)^{1/2} \,\big|\, T > t\right)$$

$$\geq \mathbf{E}_x\left(\mathbf{P}_{X_t}(T_H > s) \,\big|\, |X_t^d| \geq t^{1/4}\mathbf{P}_x(T > t)^{1/2}, T > t\right) \left(1 - \sqrt{\frac{2}{\pi}} \frac{1}{\left(\sqrt{t}\mathbf{P}_x(T > t)^{1/2}\right)}\right) \ .$$

From (3) we have

$$\mathbf{P}_x(T_H > s) = \frac{1}{\sqrt{2\pi s}} \int_0^\infty \left(e^{(x_d - y_d)^2/2s} - e^{(x_d + y_d)^2/2s}\right) dy_d$$

$$= 1 - \frac{2}{\sqrt{2\pi}} \int_{x_d/\sqrt{s}}^\infty e^{-u^2/2} du \geq 1 - \frac{2}{\sqrt{2\pi}} \frac{\sqrt{s}}{x_d} e^{-x_d^2/2s} \ .$$

Therefore

$$1 \geq \frac{\mathbf{P}_x(T > t + s)}{\mathbf{P}_x(T > t)} \geq$$

$$\left(1 - \frac{2}{\sqrt{2\pi}} \frac{\sqrt{s}}{\left(\sqrt{t}\mathbf{P}_x(T > t)\right)^{1/2}} e^{-\sqrt{t}\mathbf{P}_x(T>t)/2s}\right) \left(1 - \sqrt{\frac{2}{\pi}} \frac{1}{\left(\sqrt{t}\mathbf{P}_x(T > t)\right)^{1/2}}\right) \ ,$$

which tends to 1 when $t$ tends to infinity.

For a fixed $y \in \Omega$, using reflection arguments, Harnack inequality and monotonicity in $t$ to get away from the boundary, one can now show as in [C.M.SM.1] that as a function of $x$, any accumulation point of $\mathbf{P}_x(T > t)/\mathbf{P}_y(T > t)$ is a non negative harmonic function which is zero at the boundary and equal to one for $x = y$. The conclusion of Theorem 2 follows in that case since the cone $\mathcal{P}_\Omega$ is one dimensional.

**Appendix: a symmetry estimate.**

Let $\vec{n}$ be a unit vector in $\mathbf{R}^{d-1}$. Let $y \in \Omega$ with $y_d \neq 0$, we will assume for definiteness $y_d > 0$, the argument is the same for $y_d < 0$. Several quantities below will depend on $y$ and $\vec{n}$, but since we will not vary these two objects we will not mention this dependence explicitly. We define a domain $\Omega^+$ by

$$\Omega^+ = \left\{ x = (\vec{x}, x_d) \in \Omega \mid |x_d| \leq (\vec{x} - \vec{y}) \bullet \vec{n} \right\} .$$

We also define two symmetry operators $S^\pm$ by

$$S^+(x) = (\vec{x} + (x_d - (\vec{x} - \vec{y}) \bullet \vec{n})\vec{n}, (\vec{x} - \vec{y}) \bullet \vec{n}) ,$$

and

$$S^-(x) = (\vec{x} - (x_d + (\vec{x} - \vec{y}) \bullet \vec{n})\vec{n}, -(\vec{x} - \vec{y}) \bullet \vec{n}) ,$$

Finally, we define the function $q_t(x, y)$ by

$$q_t(x, y) = p_t(S^+(x), y) + p_t(S^-(x), y) - p_t(x, y) .$$

The following Lemma gives an estimate on $p_t(x, y)$.

**Lemma A.** *For any $x \in \Omega^+$ and any $t > 0$, we have $q_t(x, y) \geq 0$, hence*

$$p_t(x, y) \leq p_t(S^+(x), y) + p_t(S^-(x), y)$$

A similar result holds for the domain

$$\Omega^- = \left\{ x \in \Omega \mid -|x_d| \geq (\vec{x} - \vec{y}) \bullet \vec{n} \right\} .$$

The proof is based on the maximum principle. We first observe that as a function of $t$ and $x$, $q_t(x, y)$ is a solution of the heat equation. Moreover since $y_d > 0$ we have

$$q_0(x, y) = \delta(S^+(x) - y) = \delta(x - S^+(y)) .$$

The boundary of $\Omega^+$ can be decomposed in three pieces. There is first the piece corresponding to $x_d = (\vec{x} - \vec{y}) \bullet \vec{n}$, where we have $S^+(x) = x$, and therefore

$$q_t(x, y) = p_t(S^-(x), y) > 0 .$$

The second piece corresponds to $x_d = -(\vec{x} - \vec{y}) \bullet \vec{n}$, and on this hyperplane we have $S^-(x) = x$, therefore

$$q_t(x, y) = p_t(S^+(x), y) > 0 .$$

the last piece corresponds to

$$D \times \{x_d = 0\} \cap \{|x_d| \leq (\vec{x} - \vec{y}) \bullet \vec{n}\} .$$

On this set we have $p_t(x, y) = 0$, and therefore

$$q_t(x, y) = p_t(S^+(x), y) + p_t(S^-(x), y) > 0 .$$

We now conclude the proof using the Phragmén-Lindelöf version of the maximum principle (see [P.W.]).

*Acknowledgments.* The authors are grateful to the ECOS-CONYCIT grants CE99E09 and CE99E10, DIM and CMM of the Universidad de Chile, Fondap in applied Mathematics, the CNRS, University of Dijon and Ecole Polytechnique for support and hospitality.